\documentclass{article}

\usepackage{delarray,verbatim,enumerate,a4wide}
\usepackage{amsmath,amsthm,amstext,amsbsy,amssymb,amsfonts,amscd}
\usepackage[all]{xy}
\usepackage[frenchb]{babel}
\usepackage[T1]{fontenc}
\usepackage[utf8]{inputenc}
\usepackage[small]{titlesec}
\usepackage{scalerel}
\usepackage{color}
\definecolor{vert}{rgb}{0.1,0.5,0.2}
\usepackage[colorlinks,linkcolor=blue,urlcolor=blue,citecolor=vert,linktocpage]{hyperref}

\usepackage{calligra}
\DeclareFontShape{T1}{calligra}{m}{n}{<->s*[0.95]callig15}{}
\DeclareMathAlphabet{\mathscr}{T1}{calligra}{m}{n}

\newtheorem{Th}{Théorème}
\newtheorem{Lem}[Th]{Lemme}
\newtheorem{Prop}[Th]{Proposition}
\newtheorem{Cor}[Th]{Corollaire}

\newtheorem{Sco}[Th]{Scolie}
\newtheorem{Def} [Th]{Définition}

\newtheorem*{Th*}{Théorème}
\newtheorem*{ThA}{Théorème A}
\newtheorem*{ThB}{Théorème B}

\def\Preuve{\noindent {\it Preuve.~}}
\def\PreuveTh{\noindent {\it Preuve du Théorème.~}}

\def\RR{\mathbb R}		\def\QQ{\mathbb Q}	
	\def\ZZ{\mathbb Z}	
\def\F2{\mathbb{F}_2}	\def\Z2{\mathbb{Z}_2}		
\def\Zl{\mathbb{Z}_\ell} 		

 		\def\P{\mathcal  P}		\def\U{\mathcal  U}	\def\F{\mathcal  F}
\def\J{\mathcal  J}  		\def\R{\mathcal  R}	\def\D{\mathcal  D}
 	  	\def\Cl{\mathcal  C \!\ell}
\def\E{\mathcal  E}		\def\T{\mathcal  T}		\def\	

\def\q{{\mathfrak q}}	\def\p{{\mathfrak p}}		\def\a{{\mathfrak a}}
	\def\l{{\mathfrak l}}	\def\d{{\mathfrak d}}	\def\m{{\mathfrak m}}

		\def\Frob{\operatorname{Frob}}
\def\Gal{\operatorname{Gal}}		\def\Rad{\operatorname{Rad}}

\newcommand\scale[2]{\vstretch{#1}{\hstretch{#1}{#2}}}

\newcommand\si[1]{\scale{.6}{#1}}    \newcommand\su[1]{\scale{.8}{#1}}
\newcommand\ph{{\phantom{*}}}
\newcommand\ab{{\scale{.8}{\textrm{ab}}}}	\newcommand\reg{{\scale{.8}{\textrm{rég}}}}
\newcommand\aug{{\scale{.8}{\textrm{aug}}}}	\newcommand\res{{\scale{.8}{\textrm{res}}}}

\makeatletter
\newcommand*\wt[2][0.2ex]{%
        \begingroup
        \mathchoice{\wt@helper{#1}{#2}{\displaystyle}{\textfont}}
                   {\wt@helper{#1}{#2}{\textstyle}{\textfont}}
                   {\wt@helper{#1}{#2}{\scriptstyle}{\scriptfont}}
                   {\wt@helper{#1}{#2}{\scriptscriptstyle}{\scriptscriptfont}}%
        \endgroup
        #2%
}
\newcommand*\wt@helper[4]{%
        \def\currentfont{\the#41}%
        \def\currentskewchar{\char\the\skewchar\currentfont}%
        \setbox\tw@\hbox{\currentfont$#2$\currentskewchar}%
        \dimen@ii\wd\tw@
        \setbox\tw@\hbox{\currentfont$#2${}\currentskewchar}%
        \advance\dimen@ii-\wd\tw@
        \rlap{\raisebox{-#1}{$\m@th#3\kern\dimen@ii\widetilde{\phantom{#2}}$}}%
}
\makeatother

\def\%{{\scale{.8}{\infty}}}

\def\p{\mathfrak p}		\def\q{\mathfrak q}	\def\l{\mathfrak l}		

\def\k{{\boldsymbol k}}

\def\Ind{\operatorname{Ind}}	
\def\Gal{\operatorname{Gal}}	\def\Rad{\operatorname{Rad}}

\date{}

\title{\LARGE{Capitulation abélienne des groupes de classes de rayons}}

\author{Jean-François {\sc Jaulent}}

\begin{document}
\maketitle
\bigskip\bigskip

{\footnotesize
\noindent{\bf Résumé.} Nous montrons que la méthode utilisée par Bosca pour faire capituler les classes d'idéaux d'un corps de nombres s'étend aux classes de rayons dans le cas modéré. Plus précisément, nous prouvons que pour toute extension $K/\k$ de corps de nombres dans laquelle une au moins des places à l'infini se décompose complètement et tout diviseur $\m$ sans facteur carré, il existe une infinité d'extensions abéliennes $F/\k$ telles que les classes de rayons modulo $\m$ de $K$ capitulent dans $KF$. Il suit de là que les résultats de Kurihara sur la trivialité des groupes de classes des pro-extensions abéliennes maximales des corps de nombres totalement réels valent encore pour les groupes de classes de rayons dans le cas modéré.
}\smallskip

{\footnotesize
\noindent{\bf Abstract.} Building on Bosca's method, we extend to tame ray class groups the results on capitulation of ideals of a number field by composition with abelian extensions of a subfield first studied by Gras. More precisely, for every extension of number fields $K/\k$,where at least one infinite place splits completely, and every squarefree divisor $\m$ of $K$, we prove that there exist infinitely many abelian extensions $F/\k$ such that the ray class group  mod $\m$ of $K$ capitulates in $KF$. As a consequence we generalize to tame ray class groups the results of Kurihara on triviality of class groups for maximal abelian pro-extensions of totally real number fields.
}
\bigskip\bigskip

\tableofcontents
\newpage
\section*{Introduction}
\addcontentsline{toc}{section}{Introduction}

Il est reconnu que la capitulation des idéaux (i.e. le fait de devenir principal) par extension abélienne est un problème difficile. La raison en est, nous semble-t-il, que la théorie du corps de classes qui décrit précisément ces extensions abéliennes est centrée sur les questions normiques tandis que la capitulation met en jeu, elle, les propriétés du morphisme d'extension.\smallskip

 Du point de vue cohomologique, le noyau de capitulation dans une extension galoisienne est gouverné par le premier groupe de cohomologie des unités. Il est donc naturel que les propriétés des unités -- et donc le comportement des places à l'infini -- joue un rôle essentiel dans ces questions.\smallskip
 
 C'est ainsi que Gras a montré dans \cite{Gr1} que, pour une extension abélienne de corps de nombres $K/\k$ vérifiant des conditions de signatures convenables, il existe une infinité d'extensions abéliennes $F/\k$ telles que les idéaux de $K$ capitulent dans le compositum $FK$. Ce résultat a été ensuite généralisé par Bosca dans \cite{Bos} pour $K/\k$ arbitraire sous la seule condition que l'une au moins des places à l'infini se décompose complètement dans $K/\k$. Par passage à la limite inductive dans les tours abéliennes pour $\k=\QQ$, il retrouve ainsi le théorème de Kurihara \cite{Kuh} sur la trivialité du groupe de classes du compositum de $K$ avec l'extension abélienne maximale $\QQ^{ab}$ du corps des rationnels pour $K$ totalement réel et démontre en outre la conjecture de Gras (cf. \cite{Gr1}, p. 405).\smallskip
 
 La preuve de Bosca repose sur la conjonction de trois types d'arguments: la description des groupes de normes donnée par le corps de classes; le théorème de densité de Chebotarev; et des dénombrements de classes ambiges dans le cas cyclique, à la manière de Chevalley.\smallskip

Le but de la présente note est d'étendre aux groupes de classes de rayons $Cl_K^{\m_{\si{K}}}$ les résultats de Bosca en reprenant essentiellement la même stratégie -- au prix de quelques complications techniques -- dans le cas où le diviseur $\m^\ph_K$ qui définit le module de congruences est sans facteur carré, ce qui, du point de vue du corps de classes, revient à accepter de la ramification modérée.\smallskip

Le résultat que nous obtenons, en parfaite analogie avec le cas des classes d'idéaux, généralise ainsi ceux de Bosca et, partant, de Kurihara  dans ce nouveau cadre, et les redonne comme cas particuliers lorsque $\m^\ph_K$ est pris trivial.\medskip

Par exemple, pour l'extension abélienne réelle maximale de $\QQ$, le Scolie \ref{Sco} infra s'énonce ainsi:\medskip

\noindent {\bf Convention.} Étant donnés un ensemble fini $T=T_\QQ$ de places de $\QQ$ et une extension $L \subset \overline{\QQ}$, convenons de noter par $T_L$ l'ensemble des places de $L$ au-dessus de $T$ et écrivons\smallskip

\centerline{$Cl_L^{\m}=\varinjlim Cl_K^{\m}$}\smallskip

\noindent la limite inductive des groupes de classes de rayons modulo $\m^\ph_K=\prod_{\p_K^\ph\in T^\ph_K}\,\p_K^\ph$, pour $K \subset L$ de degré fini sur $\QQ$.\medskip 

Cela étant, pour tout ensemble fini $T$ de places de $\QQ$, autrement dit pour tout diviseur sans facteur carré $\m_\QQ=\prod_{p\in T}$ du corps des rationnels, il vient:

\begin{Th*}
Le sous-corps réel maximal $L=\QQ^\ab_+=\bigcup_{f>0}\QQ[\cos(2\pi/f)]$ du corps cyclotomique engendré par toutes les racines de l'unité $\QQ^\ab=\bigcup_{f>0}\QQ[\zeta^\ph_f]$ est $\m$-principal:\smallskip

\centerline{$Cl_L^\m=1$.}\smallskip

Plus généralement les extensions algébriques totalement réelles $N$ de $\QQ^\ab_+$ sont $\m$-principales.
\end{Th*}

En d'autres termes, toute classe de rayons d'un corps totalement réel modulo un diviseur sans facteur carré se trivialise par extension abélienne réelle.  La démonstration de ce résultat montre même que les extensions principalisantes $F$ peuvent être prises non-ramifiées en un ensemble fini arbitraire de places; ce qui implique l'existence d'une infinité de telles extensions non emboîtées.

\medskip

Comme expliqué en appendice, le cas sauvage, en revanche, met en jeu ultimement d'autres phénomènes qui s'opposent à la capitulation des groupes de classes infinitésimales. Par exemple, pour tout premier $\ell$ donné et sous la conjecture de Leopoldt pour $\ell$, le $\ell$-groupe des $\ell$-classes infinitésimales d'ordre fini d'un corps totalement réel ne peut capituler dans aucune de ses extensions dès lors qu'il n'est pas nul, le morphisme d'extension étant alors toujours injectif.


\newpage
\section{Énoncé des résultats et stratégie de preuve}

Rappelons qu'un corps de nombres est, par convention, une extension finie du corps $\QQ$ et qu'une place à l'infini d'un tel corps est dite se décomposer complètement dans une extension finie $K/\k$ lorsqu'elle possède {\em exactement} $[K:\k]$ prolongements à $K$ (ce qui a lieu si et seulement si ce n'est pas une place réelle dont l'un au moins des prolongements ne l'est pas).\smallskip

Le résultat principal de Bosca \cite{Bos} prouvant la conjecture de Gras peut s'énoncer comme suit:

\begin{ThA}[Bosca]\label{ThA}
Pour toute extension $K/\k$ de corps de nombres où une place à l'infini au moins se décompose complètement, il existe une extension abélienne finie $F/\k$ complètement décomposée en toutes les places à l'infini, telle que les classes de $Cl_K$ capitulent dans $L=K F$.
\end{ThA}

Sa généralisation en termes de classes de rayons s'énonce naturellement comme suit: étant donné un ensemble fini $T=T_\k$ de places d'un corps de nombres $\k$ et une extension $K$ de $\k$, convenons de noter $T_K$ l'ensemble des places de $K$ au-dessus de $T_\k$. Cela étant:

\begin{ThB}\label{ThB}
Pour toute extension $K/\k$ de corps de nombres dans laquelle une place à l'infini au moins se décompose complètement et tout ensemble fini $T$ de places non complexes de $\k$, il existe une extension abélienne  finie $F/\k$ complètement décomposée en toutes les places à l'infini, telles que les classes de rayons de $Cl_K^{\,\m}$ modulo $\m^\ph_K=\prod_{\q^\ph_K\in T^\ph_K}\q^\ph_K$  capitulent dans  $L=KF$.
\end{ThB}

Ce second résultat contient naturellement le premier qui correspond au cas particulier $T=\emptyset$. \medskip

La stratégie pour prouver le Théorème B reprend dans un contexte plus technique celle mise en œuvre par Bosca \cite{Bos} pour établir le Théorème A: comme tout groupe abélien fini est un produit direct de $\ell$-groupes cycliques, il suffit de montrer que pour tout nombre premier $\ell$ chaque classe $[\d_K^\ph]$ d'ordre $\ell$-primaire de $Cl_K^{\m_{\si{K}}}$ se trivialise dans $KF_\d$ pour une certaine $\ell$-extension abélienne $\infty$-décomposée $F_\d$ du corps de base $\k$; le compositum $F$ des $F_\d$ lorsque $\d$ parcourt un système de représentants d'une famille génératrice de classes fournissant alors un corps principalisant, toutes les classes de $Cl_K^{\m_{\si{K}}}$ capitulent par construction dans $L=KF$.\smallskip

L'idée est de représenter la classe considérée par un idéal premier convenable $\p^\ph_K$; à prendre pour $F$ une $\ell$-extension de $\k$ cyclique de degré assez grand, non-ramifiée en dehors de $\p^\ph_\k$ et non-décomposée en $\p^\ph_\k$, dans laquelle l'idéal premier $\p^\ph_\k$ au-dessous de $\p^\ph_K$ est très ramifié; puis à former le compositum $L=FK$. Par construction, l'idéal $\p^\ph_K$  est alors une grande puissance d'un idéal $\p^\ph_L$ de $L$, ambige dans l'extension cyclique $L/K$; ce qui permet de montrer que sa classe est principale, dès que l'on contrôle effectivement le nombre de classes d'ambiges dans l'extension cyclique $L/K$, ce qui se fait en imposant à $\p^\ph_K$ un certain nombre de conditions arithmétiques fortes.\smallskip

Pour cela, il est plus commode d'établir d'abord le résultat attendu dans le cas galoisien, le cas général pouvant s'en déduire ultérieurement par passage à la clôture galoisienne (Scolie \ref{Gal}).\smallskip

La situation considérée est donc la suivante:\smallskip
\begin{itemize}
\item $K/{\k}$ désigne une extension galoisienne de corps de nombres, de groupe de Galois $\Delta$, dans laquelle une au moins des places à l'infini est complètement décomposée;
\item $\ell$ est un nombre premier fixé et $T$ un ensemble fini, stable par $\Delta$, d'idéaux premiers de $K$. 
\end{itemize}\smallskip

On note $\m_K^\ph=\prod_{\q_K^\ph\in T} \q_K^\ph$ leur produit  et $E_K^{\,\m}=\{\varepsilon\in E_K^\ph \,|\, \varepsilon \equiv 1 \;{\rm mod}\; \m_K^\ph \}$ le groupe des unités de $K$ congrues à 1 modulo $\m_K^\ph$. \smallskip
En termes de représentations, l'hypothèse faite, qui revient à postuler la trivialité du sous-groupe de décomposition $\Delta_{\p_{\!\infty}}$ de l'une des places à l'infini, assure, d'après le théorème de Herbrand (cf. \cite{He1,He2} ou e.g. \cite{Mar}), que le caractère du groupe des unités\smallskip

\centerline{$\chi_{E_K}^{\phantom{l}}=\big( \sum_{\p_{\si{\infty}}} \Ind_{\Delta_{\p_{\si{\infty}}}}^\Delta 1_{\Delta_{\p_{\si{\infty}}}} \big)\;-1_\Delta^{\phantom{l}}$}\smallskip

\noindent contient le caractère d'augmentation $\chi_\Delta^\aug=\chi_\Delta^\reg\;-1_\Delta^{\phantom{l}}$. En particulier $E_K^{\,\m}$, qui est d'indice fini dans $E_K$, contient donc un sous-$\ZZ[\Delta]$-module monogène $E^{\,\varepsilon}_K=\varepsilon^{\ZZ[\Delta]}$ de caractère  $\chi_\Delta^\aug$.\smallskip

Partant d'une classe $[\d^\ph_K]$ d'ordre $\ell$-primaire de $Cl_K^{\,\m}$, on cherche un idéal premier $\p^\ph_K$ dans $[\d^\ph_K]$ et une $\ell$-extension cyclique $\p^\ph_K$-ramifiée $F$ de $\k$ convenables tels que $[\d^\ph_K]$ capitule dans $Cl_{KF}^{\,\m}$.

Le point-clé de la démonstration consiste à minorer un indice normique en imposant un comportement ad hoc de $\p^\ph_K$ dans l'extension $K_n^\varepsilon=K\big[\zeta^\ph_{\ell^n},\sqrt[\ell^{\si{n}}]{E^{\,\varepsilon}_K}\,\big]$ pour $n$ assez grand (Prop. \ref{Ambige}).

\newpage
\section{Interprétation infinitésimale des $\ell$-groupes de classes de rayons}

Soit donc $K$ un corps de nombres et $\m$ un diviseur de $K$ sans facteur carré. Écrivons\smallskip

\centerline{$K_\m^\times=\prod_{\q\mid\m} K^\times_\q$}\smallskip

\noindent le produit des groupes multiplicatifs des complétés de $K$ aux places (finies ou réelles) intervenant dans $\m$. Notons $\q_\q$ l'idéal maximal de l'anneau local des entiers de $K_\q$ et $U_\q^1=1+\q_\q$ le sous-groupe de $K_\q^\times$ formé des unités principales, si $\q$ est finie; $U_\q^1=K_\q^{\times 2}\simeq\RR_+^\times$, si $\q$ est réelle.\smallskip

Par définition (cf. e.g. \cite{Mil} \S5.1), le groupe des classes de rayons modulo $\m$ de $K$ est le quotient:\smallskip

\centerline{$Cl_K^\m=D_K^\m/P_K^\m$,}\smallskip

\noindent où $D_K^\m$ est le groupe des idéaux de $K$ étrangers à $\m$ et $P_K^\m$ le sous-groupe formé des idéaux principaux engendrés par les  $x$ de $K^\times$ qui satisfont la congruence multiplicative: $x\equiv 1 \;{\rm mod}^\times\; \m$, i.e. dont l'image canonique dans $K^\times_\m$ tombe dans le sous-groupe $U_\m^1=\prod_{\q\mid\m} U_\q^1$.\medskip

Fixons maintenant un nombre premier $\ell$; écrivons $\m=\m'\bar\m$ en mettant dans $\m'$ les places qui divisent $\ell$ et dans $\bar\m$ les autres ; puis considérons le $\ell$-sous-groupe $\,\Cl^{\,\m}_K=\Zl\otimes_\ZZ Cl^{\m}_K$ du groupe $Cl^\m_K$. Nous allons interpréter ce $\ell$-groupe en termes idéliques en nous appuyant sur le formalisme et les notations de la théorie $\ell$-adique du corps de classes tels qu'introduits dans \cite{J31}.\smallskip

Écrivons pour cela $\,\Cl^{\,\m}_K$ comme quotient du tensorisé $\ell$-adique $\D_K^{\,\m}=\Zl\otimes_\ZZ D_K^{\m}$ par son sous-module principal $\P_K^{\,\m}=\Zl\otimes_\ZZ P_K^{\m}$. 
Précisons ce dernier: pour chaque place non-complexe $\q$ de $K$, notons $\R_{K_\q}=\varprojlim \, K^\times_\q/K_\q^{\times\ell^m}$ le compactifié $\ell$-adique  du groupe multiplicatif $K^\times_{\q}$. Et notons $\,\U_{K_\q}$ son sous-groupe des unités, qui s'identifie au groupe $U_{\q}^1$ des unités principales de $K_\q$, pour $\q|\ell$; au $\ell$-sous-groupe de Sylow $\mu^\ph_{K_\q}$ du groupe des racines de l'unité contenues dans $K_\q$, pour $\q\nmid\ell$ (et donc à $\{\pm 1\}$ pour $\q$ réelle et $\ell=2$, au groupe trivial $\{1\}$ pour $\q$ réelle et $\ell\ne 2$).\smallskip

Rappelons que le tensorisé $\R_K=\Zl\otimes_\ZZ K^\times$ du groupe multiplicatif de $K$ s'injecte canoniquement dans le pro-$\ell$-groupe des idèles, défini comme le produit restreint des compactifiés $\R_{K_\q}$:\smallskip

\centerline{$\R_K\hookrightarrow \J_K=\prod^\res_\q\,\R_{K_\q}$.}\smallskip

\noindent Ainsi $\P_K^{\,\m}$ est l'image dans $\D^{\,\m}_K\simeq \prod_{\q\nmid\m}^\res\R_{K_\q}/\prod_{\q\nmid\m}^\ph\U_{K_\q} \simeq \prod_{\q\nmid\m}^\res\R_{K_\q} \prod_{\q\mid\m'}^\ph\U_{K_\q}  / \prod_{\q\nmid\bar\m}^\ph\U_{K_\q}$ de\smallskip

\centerline{$\R_K^{\,\m}=\big\{x=(x_\q)_\q\in\R_K\;|\; x_\q\in\U_\q \; \textrm{ pour } \q|\m' \textrm{ et } x_\q=1 \textrm{ pour } \q|\bar\m \big\}
= \R_K \cap  \prod_{\q\nmid\m}^\res\R_{K_\q}\prod_{\q\mid\m'}^\ph\U_{K_\q}$,}\smallskip

\noindent sous-module de $\R_K$ formé des éléments qui sont des unités aux places divisant $\m'$ et triviaux aux places divisant $\bar\m$ (on dit que que ce sont des unités en $\m'$ qui sont $\bar\m$-infinitésimales\footnote{La notion $\ell$-adique d'élément infinitésimal, introduite dans \cite{J9}, est exposée dans dans \cite{J31} et \cite{Gr2}, III.2.b.}).
Il suit:\smallskip

\centerline{$\Cl^{\,\m}_K = \D^{\,\m}_K/\P^{\,\m}_K \simeq \prod_{\q\nmid\m}^\res\R_{K_\q} \prod_{\q\mid\m'}^\ph\U_{K_\q} / \prod_{\q\nmid\m}^\ph\U_{K_\q}\,\big(\R_K\cap \prod_{\q\nmid\m}^\res\R_{K_\q}\prod_{\q\mid\m'}^\ph\U_{K_\q}\big)$.}

\noindent Puis:

\centerline{$\Cl^{\,\m}_K \simeq \prod_{\q\nmid\m}^\res\R_{K_\q} \prod_{\q\mid\m'}^\ph\U_{K_\q} \,\R_K / \prod_{\q\nmid\m}^\ph\U_{K_\q}\,\R_K$.}\smallskip

\noindent Le théorème d'approximation faible permet de remplacer au numérateur $\prod_{\q\mid\m'}^\ph\U_{K_\q}$ par $\prod_{\q\mid\m'}^\ph\R_{K_\q}$. Il vient donc finalement:\smallskip

\centerline{$\Cl^{\,\m}_K \simeq \prod_{\q\nmid\bar\m}^\res\R_{K_\q} \,\R_K / \prod_{\q\nmid\m}^\ph\U_{K_\q}\,\R_K
\simeq \prod_{\q\nmid\bar\m}^\res\R_{K_\q} / \big(\R_K \cap  \prod_{\q\nmid\m}^\res\R_{K_\q}\big) \simeq \Cl^{\,\bar\m}_K
$.}\smallskip

\begin{Prop}\label{Infinitésimaux}
Soient $\m$ un diviseur sans facteur carré d'un corps de nombres $K$, puis $\ell$ un nombre premier arbitraire et $\bar\m$ le diviseur construit sur les places de $K$ qui ne divisent pas $\ell$.\par
Le $\ell$-groupe $\,\Cl^{\,\m}_K$ des classes de rayons modulo $\m$ est indépendant de $\m'=\m/\bar\m$ et s'identifie au quotient du $\Zl$module multiplicatif $\D^{\,\bar\m}_K$ construit sur les idéaux premiers qui ne divisent par $\bar\m$ par le sous-module $\P^{\,\m}_K$ des idéaux principaux engendré par les éléments $\bar\m$-infinitésimaux $x\in\R^{\,\bar\m}_K$:\smallskip

\centerline{$\,\Cl^{\,\m}_K \simeq \,\Cl^{\,\bar\m}_K = \D^{\,\bar\m}_K/\P^{\,\bar\m}_K$.}\smallskip

Par la théorie du corps de classes, $\Cl^{\,\m}_K$ s'identifie ainsi au groupe de Galois de la $\ell$-extension abélienne $\bar\m$-ramifiée (i.e. non-ramifiée en dehors de $\bar\m$) maximale $H_K^{\bar\m}$ de $K$.
\end{Prop}

\Preuve Le corps de classes $H_K^{\m}$ associé à $\,\Cl^{\,\m}_K$ est la $\ell$-extension abélienne maximale $\m$-modérément ramifiée  (i.e. non-ramifiée en dehors de $\m$ et possiblement mais modérément aux places de $\bar\m$). Or, dans une $\ell$-extension les places divisant $\ell$ qui se ramifient le font sauvagement, les autres modérément. $H_K^{\m}=H_K^{\bar\m}$ est donc la $\ell$-extension abélienne $\bar\m$-ramifiée maximale de $K$. 

\newpage
\section{Identité des classes infinitésimales d'ambiges}

La qualification d'{\em ambiges} pour désigner les classes ou les idéaux invariants par le groupe de Galois dans une extension galoisienne de corps de nombres est traditionnelle depuis Hilbert (\cite{Hil}, \S 5.7). L'ambiguïté de cette notion tient au fait que les classes ambiges (i.e. invariantes) d'idéaux ne sont pas nécessairement représentées par des idéaux ambiges (i.e. invariants). On parle de {\em classes ambiges} dans le premier cas, de {\em classes d'ambiges} dans le second.\smallskip

Le but de la présente section est d'étendre aux $\ell$-groupes de classes infinitésimales $\,\Cl^{\,\bar\m}_L$, dans une $\ell$-extension cyclique de corps de nombres $L/K$, les calculs de Chevalley sur les classes ambiges d'idéaux dans les extensions cycliques (\cite{Chv}, pp. 402--406). Un point-clé de la démonstration est le fait que, pour $\Gamma$ cyclique, le quotient de Herbrand d'un $\ZZ[\Gamma]$-module noethérien $M$\smallskip

\centerline{$q(\Gamma,M)=\frac{|H^2(\Gamma,M)|}{|H^1(\Gamma,M)|}$}\smallskip

\noindent ne dépend que du caractère du $\QQ[\Gamma]$-module $\QQ\otimes_\ZZ M$. Pour le module des unités, il est donné par le théorème de représentation de Herbrand (cf. \cite{He1,He2} ou e.g. \cite{Mar}), ce qui conduit à la formule de Chevalley (où $d_{\q_{\si{K}}^{\si{\infty}}}(L/K)$ désigne le degré de l'extension locale attaché à la place à l'infini $\q_{\si{K}}^{\si{\infty}}$):\smallskip

\centerline{$q(\Gamma,E_L)=\frac{[L:K]}{\prod^{{}^\ph}_{\q_{\si{K}}^{\si{\infty}}}d_{\q_{\si{K}}^{\si{\infty}}}(L/K)}$}

\begin{Prop}[Classes infinitésimales d'ambiges]\label{C d'ambiges}
Soient $\ell$ un nombre premier, $L/K$ une $\ell$-extension cyclique de groupe $\Gamma$ et $\bar\m$ un diviseur $\Gamma$-invariant, étranger à $\ell$ et sans facteur carré.\par

Le nombre de $\ell$-classes de $\,\Cl^{\,\bar\m}_L$ qui sont représentées par un idéal ambige est donné par:
$$
\big(\D^{\bar\m}_L{}^\Gamma : \P^{\bar\m}_L{}^\Gamma\big)= |\,\Cl^{\bar\m}_K|\;\frac{\prod_{\q_{\si{K}}^{\si{\infty}}}d_{\q_{\si{K}}^{\si{\infty}}}(L/K)\,\prod_{\q^{\si{\circ}}_{\si{K}}\nmid{\m_{\si{K}}^\ph}} e_{\q^{\si{\circ}}_{\si{K}}}(L/K)}{[L:K]\;\big(\E^{\bar\m}_K:N_{L/K}(\E^{\bar\m}_L)\big)}
$$
\noindent où $\q_{\si{K}}^{\si{\infty}}$ parcourt les places à l'infini, $\q^{\si{\circ}}_{\si{K}}$ les places finies étrangères à $\bar\m$ ramifiées dans $L/K$; et $d_{\q^\ph_{\si{K}}}(L/K)$ et $e_{\q^\ph_{\si{K}}}(L/K)$ désignent le degré local et l'indice de ramification en $\q_{\si{K}}$.
\end{Prop}

\Preuve Reprenons les calculs effectués dans \cite{J9}, \S2 pour les classes $\ell$-infinitésimales, en nous appuyant cette fois sur la description infinitésimale des $\ell$-classes de rayons donnée dans la section précédente. Partons de l'identité: $\big(\D^{\bar\m}_L{}^\Gamma:\P^{\bar\m}_L{}^\Gamma\big) = \big(\D^{\bar\m}_L{}^\Gamma :\D^{\bar\m}_K\big)\, \big(\D^{\bar\m}_K:\P^{\bar\m}_K\big)\,/\,\big(\P^{\bar\m}_L{}^\Gamma:\P^{\bar\m}_K\big)$.\smallskip
\begin{itemize}

\item $\D^{\bar\m}_L{}^\Gamma$ est engendré par les produits $\prod_{{\q_{\si{L}}^{\si{\circ}}}|{\q_{\si{K}}^{\si{\circ}}}}{\q_{\si{L}}^{\si{\circ}}}$, où $\q_{\si{K}}^{\si{\circ}}$ décrit l'ensemble des idéaux premiers de $K$ étrangers à $\bar\m$, de sorte qu'il vient:
$(\D^{\bar\m}_L{}^\Gamma:\D^{\bar\m}_K)=\prod_{\q_{\si{K}}^{\si{\circ}}\nmid\bar\m^\ph_{\si{K}}} e_{\q_{\si{K}}^{\si{\circ}}}(L/K)$.\smallskip

\item $\D^{\bar\m}_K/\P^{\bar\m}_K$ est tout simplement le $\ell$-groupe $\,\Cl^{\bar\m}_K$.\smallskip
 
\item Reste à étudier le quotient $\P^{\bar\m}_L{}^\Gamma/\P^{\bar\m}_K$. La Proposition \ref{Infinitésimaux} nous fournit la suite exacte courte $1 \to \E^{\bar\m}_L \to \R^{\bar\m}_L \to \P^{\bar\m}_L \to 1$, où $\E^{\bar\m}_L=\E_L\cap\,\R^{\bar\m}_L$ est formé des unités de $L$ qui sont $\bar\m$-infinitésimales. Prenant la cohomologie, nous obtenons la suite exacte longue\smallskip

\centerline{$1 \to \E^{\bar\m}_L{}^\Gamma = \E^{\bar\m}_K \to \R^{\bar\m}_L{}^\Gamma = \R^{\bar\m}_K \to \P^{\bar\m}_L{}^\Gamma \to H^1(\Gamma,\E^{\bar\m}_L) \to H^1(\Gamma,\R^{\bar\m}_L) \to \cdots$}\smallskip

\noindent et finalement la suite exacte: $1 \to \P^{\bar\m}_K \to \P^{\bar\m}_L{}^\Gamma \to  H^1(\Gamma,\E^{\bar\m}_L) \to  H^1(\Gamma,\R^{\bar\m}_L)$.\smallskip

\noindent Dans celle-ci, le groupe $\R^{\bar\m}_L$ des éléments de $\R_L$ qui sont $\bar\m$-infinitésimaux est le noyau de l'épimorphisme de $\R_L$ sur $\R_{L_{\bar\m}}=\prod_{\q_|\bar\m}\R_{L_\q}$ donné par: $x=(x_\q)_\q \mapsto (x_\q)_{\q|\bar\m}$.
Partant alors de la suite  $1 \to \R^{\bar\m}_L \to \R_L \to \R_{L_{\bar\m}} \to 1$ et prenant la cohomologie, nous obtenons:\smallskip

\centerline{$1 \to \R^{\bar\m}_L{}^\Gamma = \R^{\bar\m}_K \to \R_L^\Gamma =\R_K \to \R_{L_{\bar\m}}^\Gamma = \R_{K_{\bar\m}} \to H^1(\Gamma,\R^{\bar\m}_L) \to H^1(\Gamma, \R_L) \to \cdots$.}\smallskip

\noindent de sorte que $H^1(\Gamma,\R^{\bar\m}_L)$ s'injecte dans $H^1(\Gamma, \R_L)$. Or, celui-ci est trivial, en vertu du théorème 90 de Hilbert. Il suit: $H^1(\Gamma,\R^{\bar\m}_L)=1$; puis:  $\P^{\bar\m}_L{}^\Gamma/\P^{\bar\m}_K \simeq H^1(\Gamma,\E^{\bar\m}_L)$; et finalement:\smallskip

\centerline{$\big( \P^{\bar\m}_L{}^\Gamma : \P^{\bar\m}_K\big) = |H^1(\Gamma,\E^{\bar\m}_L)|= |H^2(\Gamma,\E^{\bar\m}_L)|/q(\Gamma,\E^{\bar\m}_L)$.}\smallskip

\noindent Enfin, le quotient de Herbrand $q(\Gamma,\E^{\bar\m}_L)$ coïncide avec $q(\Gamma,\E_L)=q(\Gamma,E_L)$ puisque $\,\E_L^{\bar\m}$ est d'indice fini dans $\,\E^\ph_L$. Il vient donc:\smallskip

\centerline{$\big( \P^{\bar\m}_L{}^\Gamma : \P^{\bar\m}_K\big) = \big(\E^{\bar\m}_K : N_{L/K}(\E^{\bar\m}_L)\big) \, \frac{[L:K]}{\prod^{{}^\ph}_{\p_{\si{K}}^{\si{\infty}}}d_{\p_{\si{K}}^{\si{\infty}}}\!(L/K)}$.}
\end{itemize}

Récapitulant le tout, on obtient le résultat annoncé.

\newpage
\section{Minoration de l'indice normique des unités}

\noindent{\bf Contexte.} Pour minorer l'indice normique dans la Proposition \ref{C d'ambiges}, fixons quelques conditions:\smallskip
\begin{itemize}
\item $K/\k$ désigne désormais une extension galoisienne de corps de nombres, de groupe de Galois $\Delta$, dans laquelle une au moins des places à l'infini est complètement décomposée;
\item $\ell$ est un nombre premier fixé et $\bar\m$ un diviseur $\Delta$-stable, sans facteur carré et étranger à $\ell$;
\item $L/K=KF/K$ provient d'une $\ell$-extension cyclique $F/\k$ de groupe $\Gamma$, disjointe de $K/\k$.
\end{itemize}\smallskip

La condition de décomposition assure, d'après le théorème de Herbrand déjà cité, que le groupe $E_K^{\bar\m}=\{\varepsilon\in E_K^\ph \,|\, \varepsilon \equiv 1 \;{\rm mod} \,\bar\m_{\si{K}} \}$ contient un sous-module monogène $E^{\,\varepsilon}_K=\varepsilon^{\ZZ[\Delta]}$ de caractère  $\chi_\Delta^\aug$.\smallskip

Soit $n\ge 1$ un entier arbitraire ($n>1$, pour $\ell=2$); $\zeta_{\ell^n}$ une racine $\ell^n$-ème primitive de l'unité;  $E^{\,\circ}_{K_n}=\{\eta \in K^\times_n \,|\,   \eta^{\ell^n}\in E^{\,\varepsilon}_K \}$; et $K_n^\varepsilon=K\big[\zeta^\ph_{\ell^n},\sqrt[\ell^{\si{n}}]{E^{\,\varepsilon}_K}\,\big]$ l'extension kumérienne engendrée sur $K_n=K[\zeta^\ph_{\ell^n}]$ par les racines $\ell^n$-èmes des éléments de $E^{\,\varepsilon}_K$.

Il est bien connu que l'on a: $K^\times \cap K_n^{\times\ell^n}=K^{\times\ell^n}$ pour $\ell$ impair et $(K^\times \cap K_n^{\times 2^n}:K^{\times 2^n})\le 2$ pour $\ell=2$ (cf. e.g. \cite{Ru}, Lem. 5.7); de sorte que $\ell^{\rho_n}=(E^{\,\circ\,\ell^n}_{K_n}:E^{\,\varepsilon\,\ell^n}_K)$ est indépendant de $n$ pour $\ell$ impair, et ultimement pour $\ell=2$. Notons le $\ell^\rho$.
Pour $n\ge\rho$, l'ordre du radical $\Rad(K_n^\varepsilon/K_n)= E_K^{\,\varepsilon}/E^{\,\circ\,\ell^n}_{K_n}$ de l'extension kummérienne $K_n^\varepsilon/K_n$ est ainsi:\smallskip

\centerline{$[K_n^\varepsilon : K_n]\,=\,|\Rad(K_n^\varepsilon/K_n)|=\big( E^{\,\varepsilon}_K : E^{\,\varepsilon}_K{}^{\ell^n}\big) /(E^{\,\circ\,\ell^n}_{K_n}:E^{\,\varepsilon\,\ell^n}_K) \,=\, \ell^{([K:\k]-1)n}/\ell^{\rho}$.}\smallskip

Notons $\ell^{m}$, avec $m\ge n-\rho$, l'ordre de $\varepsilon$ dans $\,E_K^{\,\varepsilon}/E^{\,\circ\,\ell^n}_{K_n}$. Et soit $\tau_n\in\Gal(K_n^\varepsilon/K_n)$ qui vérifie:\smallskip

\centerline{$\sqrt[\ell^{\si{n}}]\varepsilon^{\,\tau_n-1}=\zeta^\ph_{\ell^m}$ pour une certaine racine $\ell^m$-ème primitive de l'unité $\zeta^\ph_{\ell^m}$.}\smallskip

Alors $\tau_n$ fixe $\sqrt[\ell^{\si{n}}]\varepsilon$ et, conjointement avec ses conjugués par $\Delta_n=\Gal(K_n/\k)$, le corps $K_n^\varepsilon$.

\begin{Prop}\label{Ambige}
Supposons que $F/\k$ est ramifiée, avec pour indice $e_{\p^\ph_\k}=\ell^n$, en une unique place $\p^\ph_\k\nmid 2\ell\bar\m^\ph_\k$ au-dessus d'un premier $p\nmid 2\ell$ qui se décompose complètement dans $\k$; que $\p_\k^\ph$ se décompose complètement dans $K_n=K[\zeta^\ph_{\ell^n}]$; et qu'il existe une place $\p_{\si{K_n}}$ de $K_n$ au-dessus de $\p^\ph_\k$ dont l'image par l'opérateur de Frobenius vérifie:\smallskip

\centerline{$\Frob(\p_{\si{K_n}}, K_n^\varepsilon/K)=\tau_n^{\ell^{n_{\si0}}}$ pour un $n_{\si{0}}\ge 0$ donné.}\smallskip

\noindent Alors le nombre de classes d'ambiges $\big(\D^{\bar\m}_L{}^\Gamma : \P^{\bar\m}_L{}^\Gamma\big)$ dans $\,\Cl_L^{\,\bar\m}$ est majoré par $\ell^{n_{\su{1}}}$ où\smallskip

\centerline{$n_{\su{1}}=h_{\si{K}}^{\su{\bar\m}}+([K:\k]-1)n_{\si{0}}+\rho$ pour $\ell\ne 2$\quad\&\quad $n_{\su{1}}=h_{\si{K}}^{\su{\bar\m}}+([K:\k]-1)n_{\si{0}}+\rho+c_{\su{L\!/\!K}}$, sinon;}\smallskip

\noindent $\ell^{\,h_{\si{K}}^{\su{\bar\m}}}=|\Cl_K^{\,\bar\m}|$ et $c_{\su{L\!/\!K}}$ est le nombre de places réelles de $K$ qui deviennent complexes dans $L$. 
\end{Prop}

\Preuve Par hypothèse les éléments de Frobenius $\Frob(\p_{\si{K_n}}^{\su{\sigma}}, K_n^\varepsilon/K_n)$ attachés aux conjuguées de la place $\p_{\si{K_n}}^\ph$ engendrent $\Gal(K_n^\varepsilon/K_n)^{\ell^{n_{\si{0}}}}$. Notons $\p_{\si{K}}^\ph$ la place de $K$ au-dessous de $\p_{\si{K_n}}^\ph$, et regardons $K$ comme plongé diagonalement dans le produit de ses complétés  $\prod_{\sigma\in\Delta}K_{\p_{\si{K}}^{\su{\sigma}}}$ au-dessus de $\p_\k^\ph$.
Du fait de l'hypothèse de décomposition,  $K_{\p_{\si{K}}^{\su{\sigma}}}$ est encore le complété de $K_n$ en la place $\p_{\si{K_n}}^{\su{\sigma}}$.
Par suite, les éléments de $K^\times$ qui sont localement des puissances $\ell^n$-èmes aux places $\p_{\si{K}}^{\su{\sigma}}$ sont globalement des puissances $\ell^{n-n_{\si{0}}}$-èmes dans $K^\times_n$. En particulier, il vient donc: $E_K^{\,\varepsilon}\cap\underset{\sigma\in\Delta}\prod K_{\p_{\si{K}}^{\su{\sigma}}}^{\times\ell^n} = E_K^{\,\varepsilon}\cap K_n^{\times\ell^{n-n_{\si{0}}}}$; et finalement:
$\big(E_K^{\,\varepsilon} : E_K^{\,\varepsilon} \cap \underset{\sigma\in\Delta}\prod U_{\p_{\si{K}}^{\su{\sigma}}}^{\ell^{\su{n}}}\big)
= [K_n[\!\sqrt[\ell^{\si{n-n_{\si{0}}}}]{E^{\,\varepsilon}_K}] : K_n] =  \ell^{([K:\k]-1)(n-n_{\si{0}})}/\ell^{\rho}$.

Maintenant, dans la formule obtenue, les groupes d'unités locales $U_{\p_{\si{K}}^{\su{\sigma}}}$ sont procycliques, puisque $p\ne 2$ étant complètement décomposé dans $K/\QQ$,  les complétés $K_{\p^{\su{\sigma}}_{\si{K}}}$ sont tous isomorphes à $\QQ_p$. 
Or, par construction, $\ell^n$ est précisément l'indice de ramification de chacune des places $\p^{\su{\sigma}}_{\si{K}}$ dans $L/K$. 
D'après la théorie locale du corps de classes,  $U_{\p_{\si{K}}^{\su{\sigma}}}^{\ell^{\su{n}}}$ est donc le sous-groupe normique  $N_{L_{\p_{\si{L}}^{\su{\sigma}}}/K_{\p_{\si{K}}^{\su{\sigma}}}}(U_{\p_{\si{L}}^{\su{\sigma}}})$ attaché à l'extension locale $L_{\p_{\si{L}}^{\su{\sigma}}}/K_{\p_{\si{K}}^{\su{\sigma}}}$. Il en résulte qu'on a:\smallskip

\centerline{$\big(E_K^{\,\varepsilon} : E_K^{\,\varepsilon} \cap \underset{\sigma\in\Delta}\prod U_{\p_{\si{K}}^{\su{\sigma}}}^{\ell^{\su{n}}}\big)= \big(E_K^{\,\varepsilon} : E_K^{\,\varepsilon} \cap \underset{\sigma\in\Delta}\prod N_{L_{\p_{\si{L}}^{\su{\sigma}}}/K_{\p_{\si{K}}^{\su{\sigma}}}}(L_{\p_{\si{L}}^{\su{\sigma}}}^\times)\big) = \big(E^{\,\epsilon}_K:E^{\,\epsilon}_K\cap N_{L/K}(L^\times)\big)$,}\smallskip

\noindent en vertu du principe de Hasse (cf. \cite{Mil}, Th. 3.1), puisque dans l'extension cyclique $L/K$ les éléments globaux qui sont  normes locales sont exactement les normes globales et que, par ailleurs, les unités sont banalement normes locales aux places non-ramifiées.
Il suit donc:\smallskip

\centerline{$\big(E_K^{\bar\m}:N_{L/K}(E_L^{\bar\m})\big) \le
 \big(E_K^{\bar\m}:E_K^{\bar\m}\cap N_{L/K}(L^\times)\big)  \le
 \big(E^{\,\epsilon}_K:E^{\,\epsilon}_K\cap N_{L/K}(L^\times)\big)=  \ell^{([K:\k]-1)(n-n_{\si{0}})}/\ell^{\rho}$.}\smallskip

\noindent Et finalement:\smallskip

\centerline{$\big(\D^{\bar\m}_L{}^\Gamma : \P^{\bar\m}_L{}^\Gamma\big)
\le |\,\Cl^{\bar\m}_K|\;\prod_{\q_{\si{K}}^{\si{\infty}}}d_{\q_{\si{K}}^{\si{\infty}}}(L/K) \;\frac{\ell^n}{[L:K]}\,\ell^{([K:\k]-1)n_{\si{0}}+\rho}
\le  |\,\Cl^{\bar\m}_K| \,2^{c_{\su{L\!/\!K}}}\,\ell^{([K:\k]-1)n_{\si{0}}+\rho}$.}

\newpage
\section{Construction de l'extension principalisante}

Récapitulons: étant donnée une extension galoisienne $K/\k$ de corps de nombres complètement décomposée en au moins une place à l'infini, un nombre premier $\ell$ et un diviseur $\m_K$ de $K$ sans facteur carré, stable par $\Delta=\Gal(K/\k)$, nous avons défini $\bar\m_K$ en écartant les places divisant $\ell$ et fait choix d'un sous-$\ZZ[\Delta]$-module $E_K^{\,\varepsilon}$ de $E_K^{\,\m}$ de caractère $\chi_\Delta^\aug$, ce qui détermine une constante $\ell^\rho$.\smallskip

Posons maintenant $\ell^{n_{\si{0}}}=[H_K\cap K_\%:K]\,|\mu^\ph_K|$, où $H_K$ désigne le $\ell$-corps de classes de Hilbert de $K$ (i.e. sa $\ell$-extension abélienne non-ramifiée maximale), $K_\%$ la $\Zl$-extension cyclotomique de $K$, et $\mu_K^\ph$ le $\ell$-sous-groupe de Sylow du groupe des racines de l'unité contenues dans $K$.
\smallskip


Étant donnés une classe $[\d_K]$ d'ordre $\ell$-primaire dans $Cl_K^\m$ et un entier $n \ge n_{\su{1}}$, nous cherchons une place $\p_\k\nmid\bar\m_\k$ de $\k$ au-dessus d'un premier $p\nmid 2\ell$ de $\QQ$ complètement décomposé dans $\k$ et une place $\p^\ph_K$ de $K$ au-dessus de $\p_\k$ qui satisfassent les quatre conditions suivantes:\smallskip
\begin{itemize}
\item[($i$)] $\;\p_\k^\ph$ est  complètement décomposée dans $K_n=K[\zeta_{\ell^n}]$;\smallskip

\item[($ii$)] l'une des places au-dessus $\p^\ph_{\su{K_n}}$ dans $K_n$ est d'image $\tau_n^{\ell^{n_{\si0}}}$ dans $\Gal\big(K_n\big[\!\sqrt[\ell^{\si{n}}]{E^{\,\varepsilon}_K}\,\big]/K_n\big)$;\smallskip

\item[($iii$)] $\p_K^\ph$ représente la classe $[\d_K]$ de $Cl_K^\m$ (autrement dit a même image dans $\Gal(H_K^{\bar\m}/K)$);\smallskip

\item[($iv$)] il existe une $\ell$-extension cyclique $F/\k$, ramifiée uniquement en $\p_\k$ avec pour indice $\ell^n$.\medskip
 \end{itemize}
 
 Examinons cette dernière condition. Par la théorie $\ell$-adique du corps de classes (cf. \cite{J31}, \S2.2), le groupe de Galois de la $\ell$-extension abélienne $\p$-ramifiée $\infty$-décomposée maximale $H_\k^\p$ du corps $\k$ relativement à sa sous-extension non-ramifiée maximale $H^\ph_\k$ est donné par l'isomorphisme:\smallskip
 
 \centerline{$\Gal(H_\k^\p/H^\ph_\k)\simeq \big(\R_\k\prod_{\q}\U_{\k_\q}\prod_{\q\mid_\infty}\R_{\k_\q}\big) /
 \big(\R_\k\prod_{\q\ne\p}\U_{\k_\q}\prod_{\q\mid_\infty}\R_{\k_\q}\big)    \simeq \mu_{\k_\p}/s_\p(\E_\k) $,}\smallskip
 
 \noindent où $\,\U_{\k_\p}=\mu_{\k_\p}$ est le $\ell$-groupe des racines de l'unité dans $\k_\p$ et $s_\p(\E_\k)$ l'image locale du groupe des unités globales.
Or, le quotient obtenu est cyclique d'ordre $\ell^d$ pour un $d\ge n$ si et seulement si le complété $\k_\p$ contient les racines $\ell^n$-ièmes de l'unité et si les éléments de $\E_\k$ sont des puissances $\ell^n$-ièmes locales dans $\k_\p$; ce qui a lieu dès que la place $\p_\k^\ph$ est complètement décomposée dans l'extension $\k \big[ \zeta_{\ell^n},\!\sqrt[\ell^n]{E_\k}\,\big]/\k$.
Le sous-groupe d'inertie $I=I(\p_\k^\ph,H_\k^\p/\k)$ est alors d'ordre $\ell^d$ pour un $d\ge n$; et il suffit de prendre pour $F$ le corps des points fixes de n'importe quel sous-groupe $B$ de $A=\Gal(H_\k^\p/\k)$ qui rencontre trivialement $I$ et définit un quotient cyclique $A/B$ pour avoir $I(\p_\k^\ph,F/\k)=BI/B\simeq I$.
En fin de compte, l'existence de $F$ est assurée si l'on remplace $(i)$ par:\smallskip
\begin{itemize}
\item[$(i')$]  $\p^\ph_\k$ se décompose complètement dans $K_n\big[\!\sqrt[\ell^{\si{n}}]{E_\k}\,\big] = K\big[\zeta_{\ell^{\su{n}}},\!\sqrt[\ell^{\si{n}}]{E_\k}\,\big]$.
\end{itemize}\smallskip

\noindent Et tout le problème est alors de s'assurer de la compatibilité des trois conditions ($i'$), ($ii$) et ($iii$).\smallskip

D'un côté, les extensions kummériennes $K_n\big[\!\sqrt[\ell^{\su{n}}]{E_\k }\,\big]/K_n$ et $K_n\big[\!\sqrt[\ell^{\su{n}}]{ E_K^{\,\varepsilon}}\,\big]/K_n$ étant linéairement disjointes, puisque  $E_\k$ et $E_K^{\,\varepsilon}$ sont en somme directe par construction, $\tau_n\in\Gal\big(K_n\big[\!\sqrt[\ell^{\si{n}}]{E^{\,\varepsilon}_K}\,\big]/K_n\big)$ se relève  dans $Gal\big(K_n\big[\!\sqrt[\ell^n]{E_\k E_K^{\,\varepsilon}}\,]/K_n\big[\!\sqrt[\ell^n]{E_\k }\,\big]\big)$ en un élément $\bar\tau_n$.\smallskip

D'un autre côté, l'intersection $H_K^E$ de $H_K^{\bar\m}$ avec $K_n\big[\!\sqrt[\ell^n]{E_\k E_K^{\,\varepsilon}}\,]=K\big[\zeta_{\ell^{\su{n}}},\!\sqrt[\ell^n]{E_\k E_K^{\,\varepsilon}}\,]$ est simultanément $\bar\m$-ramifiée (comme $H_K^{\bar\m}$) et $\ell$-ramifiée (car engendrée par des racines $\ell^n$-èmes d'unités) donc non-ramifiée (puisque $\bar\m$ et $\ell$ sont étrangers), i.e. contenue dans $H_K$. Par abélianité, il suit:\smallskip

\centerline{$H_K^E = H_K^{\bar\m} \cap K_n\big[\!\sqrt[\ell^{\su{n}}]{E_\k E_K^{\,\varepsilon}}\,] \subset H_K \cap K_n\big[\!\sqrt[\ell^{\su{m_{\si{K}}}}]{E_\k E_K^{\,\varepsilon}}\,]$, où $\ell^{m_{\su{K}}^\ph}$ mesure l'ordre du $\ell$-groupe $\mu_K^\ph$.}\smallskip

\noindent En particulier, l'exposant du $\ell$-groupe $\Gal(H_K^E/K)$ est majoré par $[H_K\cap K_\%:K]\,|\mu^\ph_K|= \ell^{n_{\si{0}}}$.\smallskip

Il suit de là que l'élément $\bar\tau_n^{\ell^{n_{\si0}}}$ regardé dans $Gal\big(K_n\big[\!\sqrt[\ell^n]{E_\k E_K^{\,\varepsilon}}\,]/K\big)$ fixe $H_K^E$. De même, si la classe $[\d]$ est une puissance $\ell^{n_{\si{0}}}$-ème dans $\,\Cl_K^{\,\bar\m}$, son image $\delta$ dans $\Gal(H_K^{\bar\m}/K)$ fixe également $H_K^E$. Ainsi $\bar\tau_n^{\ell^{n_{\si0}}}$ comme $\delta$, qui coïncident alors sur $H_K^E = H_K^{\bar\m} \cap K_n\big[\!\sqrt[\ell^{\su{n}}]{E_\k E_K^{\,\varepsilon}}\,]$, proviennent d'un même élément $\alpha$ de $\Gal\big(H_K^{\bar\m}\big[\zeta_{\ell^{\su{n}}},\!\sqrt[\ell^{\su{n}}]{E_\k E_K^{\,\varepsilon}}\,]/K\big)$.
Et le théorème de densité de Chebotarev (cf. \cite{Chb} ou  e.g. \cite{Mil}, Th. 7.4) appliqué dans la clôture galoisienne de l'extension  $H_K^\m \big[\zeta_{\ell^n}, \sqrt[\ell^n]{E_\k E_K^{\,\varepsilon}}\,\big]/\QQ$ nous assure l'existence d'une infinité de premiers $p$ possédant une place $\p_K^\ph$ au-dessus qui satisfait les conditions requises $(i')$, $(ii)$ et $(iii)$. Ainsi:

\begin{Prop}\label{Construction}
La construction de $F$ est possible dès lors que la classe $[\d_K]$ est une puissance $\ell^{n_{\si{0}}}$-ième dans $Cl_K^\m$, où $\ell^{n_{\si{0}}}$ désigne le produit du degré $[H_K\cap K_\infty:K]$ par l'ordre $\ell^{m_{\si{K}}}$ de $\mu^\ph_K$.
\end{Prop}

\newpage
\section{Preuve du résultat principal}

\begin{Th}\label{ThGalois}
Soit $K/\k$ une extension galoisienne de corps de nombres dans laquelle l'une au moins une place à l'infini se décompose complètement; $\ell$ un nombre premier; et $\bar\m_K$ un diviseur de $K$ sans facteur carré, étranger à $\ell$ et stable par $\Delta=\Gal(K/\k)$.
Pour chaque classe $[\d_K]\in Cl_K^\m$ d'ordre $\ell$-primaire du groupe de classes de rayons modulo $\m$, il existe une  $\ell$-extension abélienne $\infty$-décomposée\footnote{Quand une place réelle devient complexe par extension, certains auteurs parlent de ramification à l'infini; d'autres, d'inertie. Pour éviter toute ambiguïté, nous précisons que les places de $\k$ ne se complexifient pas dans $F$.} $F/\k$ telles que la classe $[\d_K]$ se trivialise dans le groupe de classes de rayons $Cl_L^{\bar\m}$ du compositum $L=FK$ ataché au diviseur sans facteur carré $\bar\m_{\su{L}}=\prod_{\q_{\su{L}}|\bar\m_{\su{K}}}\q_{\su{L}}$.\par
L'extension $F/\k$ peut être prise non-ramifiée en tout ensemble fini donné de places $\p_{\su{\k}}\nmid\ell$.
\end{Th}

La première étape consiste à se ramener au cas où $[\d_K]$ est une puissance $\ell^{n_{\si{0}}}$-ième dans $Cl_K^\m$.
Définissons pour cela $n_{\su{K}}$ en notant $\ell^{n_K}=|\mu^\ph_{H_K \cap K_\%}|$ (et $\ell^{n_K}=|\mu^\ph_{H_K \cap K_\%[i]}|$ dans le cas spécial $\ell=2$ et $K[i]/K$ non-ramifiée $\infty$-décomposée). Cela étant, nous avons:

\begin{Lem}\label{Puissance}
Pour tout entier $n\ge n_{\su{K}}$, il existe une $\ell$-extension abélienne $F'$ de $\QQ$ telle que le compositum $K'=KF'$ vérifie les trois conditions:\par
\centerline{$H_{K'}\cap K'_\infty=K'$, \qquad $|\mu^\ph_{K'}|= \ell^{n_{\su{K}}}\le \ell^n$, \qquad $[\d_{K'}] \in (Cl_{K'}^{\bar\m})^{\ell^n}$, où $[\d_{K'}]$ est l'étendue de $[\d_K]$.}
\end{Lem}

\Preuve Remplacer $K$ par la sous-extension non-ramifiée maximale $H_K\cap K_\infty$ de sa $\Zl$-extension cyclotomique, revient à composer $K$ avec un étage fini de la $\Zl$-extension cyclotomique $\QQ_\infty$ de $\QQ$.\par

Cela fait, le corps obtenu vérifie par construction la condition $H_{K}\cap K_\infty=K$. Par suite, le théorème de Chebotarev  appliqué dans la clôture galoisienne de l'extension $H_{K}^\m[\zeta_{2\ell^n}]$ nous assure l'existence d'un premier $\p_{\su{K}}'\nmid \ell\bar\m_{\su{K}}$ de $K$ dans $[\d_K^\ph]$ (i.e. d'image donnée dans $\Gal(H_K^{\bar\m}/K)$) au-dessus d'un $p'$ complètement décomposé dans $K[\zeta_{2\ell^n}]/\QQ$ (donc vérifiant la congruence $p' = 1 \mod \l^n$):
\begin{itemize}
\item pour $\ell$ impair, cela résulte immédiatement de la disjonction $H_K^{\bar\m} \cap K[\zeta_{2\ell^n}] =K$;
\item pour $\ell=2$, si l'extension $K[i]/K$ est non-ramifiée et $\infty$-décomposée, le résultat vaut encore à condition de remplacer $K/\k$ par l'extension $K[i]/\k$ qui vérifie alors les mêmes propriétés.
\end{itemize}
De ce fait, le sous-corps totalement réel du corps cyclotomique $\QQ[\zeta_{p'}]$ contient un unique sous-corps $F'$ qui est cyclique de degré $\ell^n$ et totalement ramifié en $p'$. Comme $p'$ est pris complètement décomposé dans $K/\QQ$, la place $\p'_{\su{K}}$ au-dessus est totalement ramifiée dans l'extension composée $KF'/K'$; de sorte que l'étendue $[\d_{KF'}]$ de $[\d_K]$  est bien une puissance $\ell^n$-ème dans $Cl_{KF'}^{\bar\m}$.\par

Enfin, $K'/K$ étant totalement ramifiée en $p'$, on a $\mu_{K'}^\ph=\mu_K^\ph$, sauf si $K[i]$ vient remplacer $K$.
\medskip

\PreuveTh D'après le Lemme \ref{Puissance}, nous pouvons supposer, sans restreindre la généralité, $H_K\cap K_\infty=K$ et $[\d_K] \in (Cl_K^{\bar\m})^{\ell^{m_{\si{K}}}}$, où $\ell^{m_{\si{K}}}$ désigne l'ordre du $\ell$-groupe $\mu^\ph_K$. Posons $n_{\su{0}}=m_{\si{K}}$ et définissons $n_{\su{1}}$ comme plus haut. La Proposition \ref{Construction} nous assure alors pour tout $n\ge n_{\su{0}}$ l'existence d'une $\ell$-extensions cyclique $\infty$-décomposée $F$ de $\k$ ramifiée en un unique premier $\p_\k$ avec pour indice $e_{\p^\ph_\k}(F/\k)=\ell^n$, qui satisfait les conditions $(i)$, $(ii)$ et $(iii)$ de la section précédente.

En particulier la classe $[\d^\ph_K]$ est représentée par l'un des $[K:\k]$ idéaux $\p^\ph_K$ au-dessus de $\p_\k$, lequel se ramifie dans l'extension composée $L/K=KF/K$ avec pour indice $\ell^n$, de sorte qu'on a: $\p^\ph_K=\a_L^{\ell^n}$ pour un idéal $\a^\ph_L$ invariant par $\Gamma=\Gal(L/K)$. L'étendue $[\d^\ph_L]=[\a^\ph_L]^{\ell^n}$ de $[\d^\ph_K]$ à $L$ est ainsi la puissance $\ell^n$-ème de la classe d'un idéal ambige. D'après la Proposition \ref{Ambige}, c'est donc la classe principale dès qu'on a $n\ge n_{\su{1}}$; et $[\d^\ph_K]$ se trivialise alors dans $Cl_L^{\bar\m}$.

\begin{Sco}\label{Gal}
La conclusion du Théorème vaut encore lorsque $K/\k$ n'est pas supposée galoisienne.
\end{Sco}

\Preuve 
Partons d'une classe d'ordre $\ell$-primaire $\delta_{\su{K}} \in Cl_K^{\bar\m}$; introduisons la clôture galoisienne $\bar K/\k$ de $K/\k$ et notons $\ell^s$ la $\ell$-partie de $[\bar K:K]$. D'après le Lemme \ref{Puissance}, quitte à remplacer $K$  par $KF'$ pour une $\ell$-extension abélienne $F'$ de $\QQ$ (ce qui remplace $\bar K$ par $\bar KF'$ avec $[\bar KF':KF'] \le [\bar K:K]=\ell^s$), nous pouvons supposer que $\delta_{\su{K}}$ est une puissance $\ell^s$-ème dans $Cl_K^{\bar\m}$, donc la norme $N_{\su{\bar K\!/\!K}}(\delta_{\su{\bar K}})$ d'une classe de $Cl_{\bar K}^{\bar\m}$.
 Donnons-nous maintenant un corps principalisant $F$ pour $\delta_{\bar K}$; et écrivons que l'étendue $j_{\su{F\!\bar K\!/\!\bar K}}(\delta_{\su{\bar K}})$ est la classe triviale $1_{\su{F\!\bar K}}$. Notant $H=FK \cap \bar K$ et $\delta_{\su{H}}=N_{\su{\bar K\!/\!H}}(\delta_{\su{\bar K}})$, nous avons:\smallskip

\centerline{$j_{\su{F\!K\!/\!H}}(\delta_{\su{H}}) = j_{\su{F\!K\!/\!H}} (N_{\su{\bar K\!/\!H}}(\delta_{\su{\bar K}})) = N_{\su{F\!\bar K\!/\!F\!K}} (  j_{\su{F\!\bar K\!/\!\bar K}}(\delta_{\su{\bar K}})) = N_{\su{F\!\bar K\!/\!F\!K}}(1_{\su{F\!\bar K}}) =  1_{\su{F\!K}}$;}\smallskip

\noindent et $j_{\su{H\!/\!K}}(\delta_{\su{K}}) = j_{\su{H\!/\!K}} (N_{\su{H\!/\!K}}(\delta_{\su{H}})) = \prod_{\sigma\in\Gal(H\!/\!K)}\delta_{\su{H}}^\sigma$ capitule donc dans $FK$.
\newpage
\section{Conséquences arithmétiques}

\begin{Def}
Soient $\k$ un corps de nombres et $T=T_\k$ un ensemble fini de places de $\k$. Pour chaque extension algébrique $K$ de $\k$ (non nécessairement de degré fini sur $\k$), convenons de noter $T_K$ l'ensemble des places de $K$ au-dessus de $T$. 
Nous disons que $K$ est $T_K^\circ$-principal lorsque son groupe des classes $T_K^\circ$-infinitésimales 
 (défini comme limite inductive des groupes de classes de rayons $Cl_H^\m$ modulo $\m_{\su{H}}=\prod_{\q_{\su{H}}\in T_H}\q_{\su{H}}$ associés aux sous-extensions $H$ de degré fini) est trivial:
 
\centerline{$Cl_K^{T_K^\circ}= \varinjlim Cl_H^\m =1$}
\end{Def}

Cela posé, appliquons d'abord le Théorème principal (Th. B) avec $\k=\QQ$.

\begin{Sco}\label{Sco}
Pour tout corps de nombres totalement réel $K$ et tout ensemble fini $T=T_\QQ$ de places de $\QQ$, il existe une infinité d'extensions abéliennes réelles $F/\QQ$ telles que les classes de rayons modulo $\m_K^\ph=\prod_{\q_K^\ph\in T_K^\ph}\q_K^\ph$ de $Cl_K^{\m_{\si{K}}}$ capitulent dans  $Cl_{F\!K}^{\m_{\si{F\!K}}}$ (avec $\m_{F\!K}^\ph=\prod_{\q_{F\!K}^\ph\in T_{F\!K}^\ph}\q_{F\!K}^\ph$).\par
En d'autres termes, le groupe $Cl_K^{T_K^\circ}$ des classes $T_K^\circ$-infinitésimales de $K$ capitule dans le sur-corps $K[\cos(2\pi/n)]$ pour une infinité de $n$.
\end{Sco}

Passant à la limite inductive, nous obtenons la généralisation suivante, en termes de classes infinitésimales modérées, du Théorème principal de Kurihara (cf. \cite{Kuh}, Th. 1.1 \& \cite {Bos}, Cor. 1):

\begin{Cor}
Soit $\QQ^\ab_+=\bigcup_{n>0}\QQ[\cos(2\pi/n)]$ la plus grande extension abélienne réelle de $\QQ$ (i.e. le sous-corps réel maximal du corps cyclotomique $\QQ^\ab=\bigcup_{n>0}\QQ[\zeta_n]$).
Alors, pour tout ensemble fini $T_\QQ$ de places de $\QQ$, les extensions algébriques totalement réelles $K$ de $\QQ^\ab_+$ sont $T_K^\circ$-principales.
\end{Cor}




Regardons maintenant le cas relatif en distinguant suivant la signature du corps de base $\k$:

\begin{Sco}
Soit $\k$ un corps de nombres totalement réel et $T_\k$ un ensemble fini de places de $\k$.
Alors, pour tout corps de nombres $K\supset\k$ où l'une au moins des places réelles se décompose complètement, il existe une infinité d'extensions abéliennes $F/\k$ complètement décomposées en toutes les places à l'infini et telles que les classes de rayons 
de $Cl_K^{\m_{\su{K}}}$ capitulent dans $Cl_{F\!K}^{\m_{\su{F\!K}}}$.
\end{Sco}

Passant à la limite, nous obtenons la généralisation en termes de classes infinitésimales modérées d'un résultat conjecturé par Gras et démontré par Bosca (cf. \cite{Gr1}, Conj. 0.5 \& \cite{Bos}, Cor. 4):

\begin{Cor}
Soit $\k$ un corps de nombres totalement réel, $\k^\ab_+$ sa plus grande extension abélienne totalement réelle et $T_\k$ un ensemble fini de places de $\k$. Alors toute extension algébrique $K$ de $k^\ab_+$ dont la clôture galoisienne possède au moins un plongement réel est $T_K^\circ$-principale.
\end{Cor}

\Preuve D'après le scolie, pour tout $\alpha$ dans $K$, le groupe des classes $T^\circ_{\k[\alpha]}$-infinitésimales du corps $\k[\alpha]$ capitule dans le sous-corps $\k^\ab_+[\alpha]$ de $K$. Et $K$ est donc bien $T_K^\circ$-principal, comme annoncé.

\begin{Sco}
Soit $\k$ un corps de nombres qui possède au moins une place complexe 
et $T_\k$ un ensemble fini de places de $\k$.
Alors, pour tout corps de nombres $K$ qui contient $\k$, il existe une infinité d'extensions abéliennes $F/\k$ complètement décomposées en toutes les places à l'infini et telles que les classes de rayons modulo $\m_K^\ph=\prod_{\q_K^\ph\in T_K^\ph}\q_K^\ph$ de $Cl_K^\m$ capitulent dans $Cl_{F\!K}^{\m_{\si{F\!K}}}$.
\end{Sco}

Passant à la limite inductive, nous obtenons, toujours en termes de classes infinitésimales modérées, une généralisation d'un second résultat de Kurihara (cf. \cite{Kuh}, Th. A.1 \& \cite{Bos}, Cor. 3):




\begin{Cor}
Soient $\k$ un corps de nombres qui possède au moins une place complexe, $\k^\ab_+$ sa plus grande extension abélienne de $\k$ complètement décomposées en toutes les places à l'infini, et $T_\k$ un ensemble fini de places de $\k$.

Alors toute extension algébrique $K$ de $\k^\ab_+$ est $T_K^\circ$-principale: $Cl_K^{T_K^\circ}=1$.
\end{Cor}

\Preuve D'après le scolie, pour tout $\alpha$ dans $K$, le groupe des classes $T^\circ_{\k[\alpha]}$-infinitésimales du corps $\k[\alpha]$ capitule dans le sous-corps $\k^\ab_+[\alpha]$ de $K$. Et $K$ est donc bien $T_K^\circ$-principal, comme annoncé.

\newpage
\section*{Appendice}
\addcontentsline{toc}{section}{Appendice}
\medskip

 Nous nous sommes limités dans cet article au cas de la ramification modérée. On peut naturellement se demander ce qu'il en est lorsqu'on autorise de la ramification sauvage, notamment en considérant les groupes de classes $T$-infinitésimales plutôt que $T^\circ$-infinitésimales. \medskip

Fixons pour cela un nombre premier $\ell$. Une première difficulté est que le $\ell$-groupe des classes $T$-infinitésimales $\,\Cl_K^T$ n'est plus nécessairement fini dès lors que l'ensemble $T$ contient une ou plusieurs des places au-dessus de $\ell$: par exemple, si $T$ les contient toutes, la conjecture de Leopoldt postule précisément que  $\,\Cl_K^T$ est fini dans le cas où le corps $K$ est totalement réel; il est toujours infini sinon (cf. e.g. \cite{Gr2}, III \S3, ou \cite{J9,J31}). Plus généralement, si la $\ell$-partie $T'$ de $T$ contient certaines des places au dessus de $\ell$ (sans nécessairement les contenir toutes), la dimension du groupe $\,\Cl_K^T$ modulo sa $\Zl$-torsion dépend de propriétés fines d'indépendance $\ell$-adique des  $\ell$-unités.\medskip

On peut bien sûr tenter de contourner cette obstruction en se concentrant sur le sous-groupe de torsion $\,\T_K^T$ de $\,\Cl_K^T$, qui est seul susceptible de capituler par $\ell$-extension. Mais, même dans ce cadre restreint, le résultat de capitulation peut être en défaut: par exemple, lorsque $T$ est exactement l'ensemble des places au-dessus de $\ell$, le morphisme d'extension $\T_K^T \rightarrow T_L^T$ est toujours injectif sous la conjecture de Leopoldt (cf. e.g. \cite{Gr2}, IV Th.2.1, ou \cite{J9,J31}), quelle que soit l'extension $L/K$. La raison profonde est que le groupe d'unités qui contrôle la capitulation est trivial dans ce dernier cas. Un phénomène analogue (mais plus compliqué) se produit pour le sous-groupe de torsion du module de Bertrandias-Payan $\,\T^{bp}_K$ (cf. \cite{Gr5,J56,Ng2}), pour lequel le groupe en question se réduit aux racines de l'unité.


\addcontentsline{toc}{section}{Bibliographie}
\def\refname{\normalsize{\sc  Références}}
\bigskip
{\small
\noindent{\sc Adresse:}
Univ. Bordeaux \& CNRS,
Institut de Mathématiques de Bordeaux, UMR 5251,
351 Cours de la Libération,
F-33405 Talence cedex

\noindent{\sc Courriel:}
{\tt jean-francois.jaulent@math.u-bordeaux.fr}
 
 \noindent{\sc url:} \url{https://www.math.u-bordeaux.fr/~jjaulent/}
}


{\footnotesize
\begin{thebibliography}{tt}

\bibitem{Bos} {\sc S. Bosca},
\textit{Principalization of ideals in Abelian extensions of number fields},
Int. J. Number Th. {\bf 5} (2009), 527--539.

\bibitem{Chb} {\sc N. Chebotarev},
\textit{Die Bestimmung der Dichtigkeit einer Menge von Primzahlen, welche zu einer gegebenen Substitutionsklasse gehören}, 
Math. Ann. {\bf 95} (1926), 191--228.

\bibitem{Chv}{\sc C. Chevalley}
\textit{Sur la théorie du corps de classes dans les corps finis et les corps locaux}
J. fac. Sci. Tokyo 2 (1933), 365--476.
 
 \bibitem{Gr1}{\sc G. Gras}, 
 \textit{Principalisation d'idéaux par extensions absolument abéliennes}
 J. Number Th. {\bf 62} (1997), 403--421.
 
\bibitem{Gr2}{\sc G. Gras}, 
\textit{Class Field Theory: From Theory To Practice},
Springer, 2005.

\bibitem{Gr5}  {\sc G. Gras},
{\em Sur le module de Bertrandias--Payan dans une $p$-extension galoisienne -- Noyau de capitulation},
Pub. Math. Besançon  2016, 25--44.

\bibitem{Gr6}  {\sc G. Gras},
{\em Invariant generalized ideal classes – structure theorems for $p$-class groups in $p$-extensions},
Proc. Math. Sc.  {\bf 127} (2017), 1--34.

\bibitem{He1} {\sc J. Herbrand},
\textit{Nouvelle démonstration et généralisation d'un théorème de Minkowski},
C.R.A.S. {\bf 191} (1930), 1282--1285.

\bibitem{He2} {\sc J. Herbrand},
\textit{Sur les unités d'un corps algébrique},
C.R.A.S. {\bf 192} (1931), 24--27.

\bibitem{Hil} {\sc D. Hilbert}
{\textit Théorie des corps de nombres algébriques},
Ann. Fac. Sci. Toulouse {\bf 2} (1910), 225--456.

\bibitem{J9} {\sc J.-F. Jaulent},
\textit{$\mathscr S\,$-classes infinitésimales d'un corps de nombres algébriques}, 
 Ann. Sci. Inst. Fourier {\bf 34} (1984), 1--27.
 

\bibitem{J31}{\sc J.-F. Jaulent},
\textit{Théorie $\ell$-adique globale du corps de classes},
J. Théor. Nombres Bordeaux {\bf 10} (1998), 355--397.

\bibitem{J56}{\sc  J.-F. Jaulent},
\textit{Sur la capitulation pour le module de Bertrandias-Payan},
Pub. Math. Besançon  2016, 45--58.

\bibitem{Kuh} {\sc M. Kurihara},
\textit{On the ideal class group of the maximal real subfields of number fields with all roots of unity},
J. European Math. Soc.{\bf 1} (1999), 35--49.

\bibitem{Mil}{\sc J.S. Milne}, 
\textit{Class field theory}, v 4.02 (2013).
\url{https://www.jmilne.org/math/CourseNotes/CFT.pdf}

\bibitem{Mar}{\sc J. Martinet}, 
\textit{Le théorème de Herbrand sur les unités},
Sém. Th. Nombres Bordeaux (1968-1968), exp. \no 10, 1--11. \url{http://www.numdam.org/article/STNB_1968-1969____A10_0.pdf}

\bibitem{Ng2} {\sc T. Nguyen Quang Do},
{\it Descente galoisienne et capitulation entre modules de Bertrandias--Payan},
Pub. Math. Besançon  2016, 59--79.

\bibitem{Ru} {\sc K.Rubin},
{\em Global units and ideal class groups},
Invent. Math. {\bf 89} (1987), 511--526.

\end{thebibliography}
}
\end{document}